\documentclass[11pt]{amsart}

\usepackage{amssymb,latexsym,euscript,amscd,amsmath}

\input xypic

\DeclareMathAlphabet{\E}{U}{eus}{m}{n}

\theoremstyle{plain}

\newtheorem{thm}{Theorem}
\newtheorem{lem}{Lemma}
\newtheorem{prop}{Proposition}

\theoremstyle{remark}

\theoremstyle{definition}

\title{A new realization of the cohomology of Springer fibers}

\author{Shrawan Kumar \and  Jesper Funch Thomsen}

\address
{Department of Mathematics, University of North Carolina, Chapel
Hill, NC 27599-3250, USA, and
Institut for Matematisk Fag, Aarhus
Universitet, Ny Munkegade, DK-8000 \AA rhus C, Denmark.}

\email{shrawan$@$email.unc.edu, funch$@$imf.au.dk}

\dedicatory{ Dedicated to Professor M.S. Raghunathan}

\begin{document}

\maketitle

\section{Introduction}
Fix a positive integer $n$ and consider the algebraic group
$G=SL_n({\mathbb{C}})$ with its Lie algebra $sl_n({\mathbb{C}})$.
For any partition $\sigma$ of $n$, let $ X_\sigma \subset G/B$ be the associated 
Springer fiber, where $B$ is the standard Borel subgroup consisting of the upper triangular matrices. By the poineering work of Springer, its cohomology $H^*(X_\sigma)$  with complex coefficients admits an action of the Weyl group
$S_n$. Subsequently, the $S_n$-algebra $H^*(X_\sigma)$ played a fundamental role in several diverse problems. The aim of this short note is to give a geometric realization of $H^*(X_\sigma)$. More specifically, we prove the following
theorem which is the main result of this note.
\vskip2ex

{\bf Theorem 2.} {\it The coordinate ring  $\mathbb{C}[(N_G(T)\cdot 
\mathcal N_{\sigma^\vee}) \cap
\mathfrak{h} ]$ of the scheme  theoretic intersection of 
$N_G(T)\cdot 
\mathcal N_{\sigma^\vee}$ with the Cartan subalgebra 
$\mathfrak{h}$ in $sl_n({\mathbb{C}})$ is isomorphic to ${\rm H}^*(X_\sigma)$ as a graded
$S_n$-algebra, where $\sigma^\vee$ is the dual partition of $\sigma$, $T$ is the maximal torus consisting of the diagonal matrices with $\mathfrak{h}:=$ Lie 
$T$,
$N_G(T)$ is its normalizer in $G$, and 
$\mathcal N_{\sigma^\vee}$ is the full nilpotent cone of the Levi component
of the parabolic subalgebra of $sl_n({\mathbb{C}})$ associated to the partition $\sigma^\vee$. }

\vskip2ex

This theorem should be contrasted with the following theorem of 
de Concini-Procesi.

\vskip2ex

{\bf Theorem}  (\cite{ConProc}, Theorem 4.4). {\it The cohomology algebra
 ${\rm H}^*(X_\sigma)$ is isomorphic, as a graded $S_n$-algebra,  with
the coordinate ring $\mathbb{C}[\overline{G \cdot M_{\sigma^\vee}}
\cap \mathfrak{h} ]$ of the scheme theoretic intersection of
$\mathfrak h$ with the closure of the $G$-orbit of
$M_{\sigma^\vee}$, where  $M_{\sigma^\vee}$ is a nilpotent matrix associated to the partition $\sigma^\vee$.}

\vskip2ex

The proof of our theorem is based on a certain characterization of the 
$S_n$-algebra ${\rm H}^*(X_\sigma)$ given in Proposition \ref{universal property}, which seems to be of independent interest. The proof of this proposition is based on some works of Bergeron-Garsia, Garsia-Haiman and Garsia-Procesi revolving around the so called $n!$ conjecture. 

Finally, it should be mentioned that the  direct analogue of our theorem
(and also the above theorem of de Concini-Procesi) for other groups 
does not hold in general. However a partial generalization of the result of 
 de Concini-Procesi is obtained by Carrell \cite{carrell}. 

\section{Notation and Preliminaries}
Fix a positive integer $n$ and consider the algebraic group
$G=SL_n({\mathbb{C}})$. By $\mathcal{N}$ we denote the full nilpotent
cone inside the Lie algebra $sl_n({\mathbb{C}})$ of $G$. The group
$G$ acts on $\mathcal{N}$ by the adjoint action with finitely many
orbits. An orbit is determined uniquely by the sizes of the Jordan
blocks of any  element in the orbit, and this sets up a one to
one correspondence between the partitions of $n$ and the $G$-conjugacy classes
 inside $\mathcal{N}$. For each partition
$\sigma : \sigma_0 \geq \sigma_1 \geq \dots \geq \sigma_m
>0$ of $n$,  we let $M_\sigma$ denote the nilpotent matrix in the Jordan normal form
with blocks of sizes $\sigma_0, \sigma_1, \dots, \sigma_m$ along the diagonal
in the stated order and starting from the upper left corner.

Let $B$ denote the Borel subgroup of $G$ consisting of the upper triangular matrices  and
let $T$ denote the group of diagonal matrices in $G$. The Lie
algebras of $B$ and $T$ will be denoted by $\mathfrak b$ and
$\mathfrak h$ respectively. For any partition $\sigma$ of $n$ we let $X_\sigma$
denote the closed subset (called the {\it Springer fiber})
$$ X_\sigma := \{ gB \in G/B : {\rm Ad}(g^{-1})M_\sigma \in \mathfrak b \}$$
of $G/B$. This can also be identified with the set of Borel subalgebras of
$sl_n({\mathbb{C}})$ containing $M_\sigma$ or with certain fibers of the
Springer resolution of the nilpotent cone.

The singular cohomology ring ${\rm H}^*(X_\sigma) ={\rm H}^*(X_\sigma, \mathbb{C})$ with complex coefficients has
an action of the symmetric group $S_n$ on $n$-letters, 
the well known Springer representation. It is known that  ${\rm
H}^*(X_\sigma, \mathbb{C})$ is zero in odd degrees, so in the
following we will consider it as a (commutative) graded algebra under rescaled
grading by assigning degree $i$ to the elements of degree $2i$.  By
the Springer correspondence, the top degree part
${\rm H}^{d_\sigma}(X_\sigma)$ is an irreducible $S_n$-module.

For the partition  $\sigma_o : 1 \geq 1 \geq \dots \geq 1$ of $n$, the variety $X_{\sigma_o}$ coincides
with $G/B$. Thus, in this case, one may  $S_n$-equivariantly identify ${\rm
H}^*(X_{\sigma_o})$   with the covariant ring ${\mathbb{C}}[Z_1, \dots, Z_n]/I$,
where $I$ is the ideal generated by the elementary symmetric functions in the
variable $Z_1, \dots,Z_n$. For a general partition $\sigma$ of $n$, the natural map
$${\rm H}^*(G/B) \rightarrow {\rm H}^*(X_\sigma)$$
is a surjective $S_n$-equivariant map \cite{spaltenstein}. This also 
follows from the result of de Concini-Procesi mentioned in the 
introduction.

\subsection{The algebra $A_\sigma$}

For any partition  $\sigma : \sigma_0 \geq \sigma_1 \geq \cdots \geq \sigma_m >0$  of $n$, let  $D_\sigma$  be the set of pairs of nonnegative
integers $(i,j)$ satisfying $i < \sigma_j$. Then $D_\sigma$ consists of $n$
elements and we fix an ordering $\{(i_s,j_s) \}_{s=1,2,\dots,n}$ of these. 
Define the polynomial
$$\Delta_\sigma = {\rm det}[X_s^{i_t}Y_s^{j_t}]_{1 \leq s,t \leq n}\in R_n,$$
where $R_n$ is the polynomial ring $ \mathbb{C}[X_1,\dots, X_n, Y_1, \dots, Y_n]$.
Observe  that, up to a sign, $\Delta_\sigma$ does not depend on the choice of
the ordering
of the elements in $D_\sigma$.

The group $S_n$ acts on $R_n$ by acting in the natural way on the two sets of
variables $X_1, \dots, X_n$ and $Y_1, \dots,Y_n$ diagonally. We think of $R_n$ as 
a bigraded $S_n$-module by counting the degrees in the two sets of variables
separately. Define a bigraded $S_n$-equivariant ideal
in $R_n$ by :
$$ K_\sigma = \{ f \in R_n : f
(\frac{\partial}{\partial X_1},\dots, \frac{\partial}{\partial X_n},
\frac{\partial}{\partial Y_1}, \dots, \frac{\partial}{\partial Y_n})
\Delta_\sigma = 0\},$$ 
where $\frac{\partial}{\partial X_i}$ and
$\frac{\partial}{\partial Y_j}$ are the usual differential operators on $R_n$.

Define now a bigraded $S_n$-algebra  $A_\sigma := R_n /
K_\sigma$. This algebra was introduced by Garsia-Haiman
 who conjectured that $A_\sigma$ has dimension $n!$ for any partition
$\sigma$ of $n$ \cite{Garsia-Haiman}.
 This conjecture, which was called  the $n!$-conjecture,  is now 
proved by  Haiman \cite{haiman}.

Let $d_1$ (resp. $d_2$) be the $X$-degree (resp. $Y$-degree) of
$\Delta_\sigma$. Then, it is easy to see that, the bigraded component
$A_\sigma^{(d_1,d_2)}$ is one
dimensional and,  moreover, if $A_\sigma^{(e_1,e_2)} \neq 0$, then 
$e_1\leq d_1, e_2\leq d_2$. Clearly,
\begin{equation}
\label{formula}  d_1 = \sum_{(i,j) \in D_\sigma}i,\,\,\,d_2 = \sum_{(i,j) \in D_\sigma}j = \sum_{s=0}^\ell
\binom{\sigma^\vee_s}{2},
\end{equation}
where $\sigma^\vee : \sigma^\vee_0 \geq \sigma^\vee_1 \geq \dots \geq \sigma^\vee_\ell
>0$ is the dual partition. 

By definition, it follows  easily that $A_\sigma$ is Gorenstein (see,
e.g., \cite{Eisenbud}, Exercise 21.7), and hence that $A_\sigma$
has a unique minimal nonzero ideal $A_\sigma^{(d_1,d_2)}$.

As explained in \cite{Garsia-Haiman}, Section 3.1, the following theorem
 follows from the results in \cite{Bergeron-Garsia}
and \cite{Garsia-Procesi}.

\begin{thm}
\label{Haiman-Garsia algebra} The subalgebra of $A_\sigma$ generated
by the images of $X_1,\dots,X_n$ is $S_n$-equivariantly isomorphic
to ${\rm H}^*(X_{\sigma^\vee})$. Similarly, the subalgebra of
$A_\sigma$ generated by the images of $Y_1, \dots, Y_n$ is
$S_n$-equivariantly isomorphic to ${\rm H}^*(X_{\sigma})$.
	If we assign degree $1$ to all the elements $X_1,\dots,X_n,
Y_1, \dots, Y_n$, then both of these isomorphisms are graded algebra 
isomorphisms. 
\end{thm}

\section{A geometric realization  of ${\rm H}^*(X_\sigma)$}

In this section we give a new geometric realization of ${\rm H}^*(X_{\sigma})$.
Recall that the socle of a ring is defined to be the sum of all its
minimal nonzero ideals. Then, with the notation from the previous
section, we have the following: 

\begin{lem}
\label{socle} The top degree $d_\sigma$ of ${\rm H}^*(X_\sigma)$
is equal to $d_2$. Moreover, the socle of ${\rm H}^*(X_{\sigma})$
coincides with the top degree part ${\rm H}^{d_2}(X_\sigma)$.
\end{lem}
\begin{proof}
Let $z$ denote any nonzero homogeneous element in ${\rm
H}^*(X_\sigma)$. By Theorem \ref{Haiman-Garsia algebra} we may
regard $z$ as the image in $A_\sigma$ of a homogeneous polynomial
$f$ in the variables $Y_1, \dots, Y_n$. As $A_\sigma$ is
Gorenstein, we can find a monomial
$$g = X_1^{\alpha_1} \cdots X_n^{\alpha_n} Y_1^{\beta_1} \cdots
Y_n^{\beta_n}$$ 
such that the image of $f \cdot g$ in $A_\sigma$
is nonzero and has the maximal degree, i.e.,  has bidegree
$(d_1,d_2)$. But then the image of $f' := f Y_1^{\beta_1} \cdots
Y_n^{\beta_n}$ in $A_\sigma$ is nonzero and, of course,  lies in the subalgebra
generated by the images of $Y_1, \dots, Y_n$. In particular, by
Theorem \ref{Haiman-Garsia algebra}, the image of $f'$ corresponds
to a nonzero element in ${\rm H}^{d_2}(X_\sigma)$  which
equals the product $z \cdot z'$ for some $z'$ in ${\rm
H}^*(X_\sigma)$.

This proves that $d_2$ equals the top degree of ${\rm
H}^*(X_\sigma)$ and that any nonzero element  of ${\rm
H}^*(X_\sigma)$ can be multiplied by an element  of ${\rm
H}^*(X_\sigma)$ to produce
 a nonzero element in the top degree. This immediately implies
the desired result.
\end{proof}

The above lemma provides us with the following characterization of the algebra
${\rm H}^*(X_\sigma)$.

\begin{prop}
\label{universal property} Let $K$ be a graded algebra with an action
of $S_n$ such that there exists a surjective $S_n$-equivariant
graded algebra homomorphism  $\phi : {\rm H}^*(X_\sigma) \rightarrow K$. Assume
further that the top
degree of $K$ is $d_\sigma$. Then $\phi$ is an isomorphism.
\end{prop}
\begin{proof}
Assume that $\phi$ is not injective. Then ${\rm ker}(\phi)$  will
meet the socle of ${\rm H}^*(X_\sigma)$ nontrivially. Thus, by Lemma \ref{socle}, the degree $d_\sigma$ part ${\rm
ker}^{d_\sigma}(\phi)$ of ${\rm ker}(\phi)$ is nonzero. But ${\rm
ker}^{d_\sigma}(\phi)$ is a submodule of the irreducible
$S_n$-module ${\rm H}^{d_\sigma}(X_\sigma)$ and hence
${\rm ker}^{d_\sigma}(\phi)$ coincides with ${\rm
H}^{d_\sigma}(X_\sigma)$. This is a contradiction, since  the top degree of 
$K$ is
 equal to $d_\sigma$ by assumption.
\end{proof}

Let $Z_T(M_\sigma)$ denote the centralizer of $M_\sigma$ in $T$,
and let $L_\sigma$ denote the centralizer in $G$ of the group
$Z_T(M_\sigma)$. In other words, $L_\sigma$ is the set of block
diagonal matrices in $G$ with blocks of sizes $\sigma_0,\sigma_1,
\dots, \sigma_m$. Then $L_\sigma$ is a reductive group with a
Borel subgroup $B_\sigma = B \cap L_\sigma$. Moreover, $M_\sigma$
is a regular nilpotent element in the Lie algebra of $L_\sigma$
and consequently the nilpotent cone $\mathcal N_\sigma$ of the Lie
algebra of $L_\sigma$ coincides with the closure
$\overline{L_\sigma  M_\sigma}$ of the $L_\sigma$-orbit of
$M_\sigma$ under the adjoint action. In the following, $N_G(T)$
will denote the normalizer of $T$ in $G$.

\begin{lem}
\label{S_n-structure} The coordinate ring $\mathbb{C}[(N_G(T)\cdot
\mathcal N_\sigma) \cap \mathfrak{h} ]$ of the scheme theoretic
intersection of the $N_G(T)$-orbit of $\mathcal N_\sigma$ with
$\mathfrak h$ is a graded $S_n$-algebra.
\end{lem}
\begin{proof}
The multiplication action of $\mathbb{C}^*$ on the Lie algebra
$sl_n(\mathbb{C})$ keeps $N_G(T)\cdot \mathcal N_\sigma$ and $\mathfrak
h$ stable. This defines the desired grading on $\mathbb{C}[(N_G(T)\cdot
\mathcal N_\sigma) \cap \mathfrak{h} ]$. The natural actions of
$N_G(T)$ on $N_G(T)\cdot \mathcal N_\sigma$ and $\mathfrak{h}$ define
an action of $N_G(T)$ on the coordinate ring $\mathbb{C}[(N_G(T)\cdot
\mathcal N_\sigma) \cap \mathfrak{h}]$. As $T$ acts trivially on
$\mathfrak h$, this gives the desired $S_n$-algebra structure by
identifying $S_n$ with the Weyl group $N_G(T)/T$.
\end{proof}

\begin{thm}
The algebra $\mathbb{C}[(N_G(T)\cdot \mathcal N_{\sigma^\vee}) \cap
\mathfrak{h} ]$ is isomorphic to ${\rm H}^*(X_\sigma)$ as a graded
$S_n$-algebra.
\end{thm}
\begin{proof}
By the result of de Concini-Procesi mentioned in the introduction, we
may identify ${\rm H}^*(X_\sigma)$ as a graded $S_n$-algebra with
the coordinate ring $\mathbb{C}[\overline{G \cdot M_{\sigma^\vee}}
\cap \mathfrak{h} ]$ of the scheme theoretic intersection of
$\mathfrak h$ with the closure of the $G$-orbit of
$M_{\sigma^\vee}$. The graded $S_n$-algebra structure on the latter
algebra is defined similarly to the graded $S_n$-algebra structure on
$\mathbb{C}[(N_G(T)\cdot \mathcal N_\sigma) \cap \mathfrak{h} ]$ (as
defined in the proof of Lemma \ref{S_n-structure}). As $N_G(T)\cdot
\mathcal N_{\sigma^\vee} = N_G(T)\cdot (\overline{L_{\sigma^\vee} \cdot 
M_{\sigma^\vee} })$
is a closed subscheme of $\overline{G \cdot M_{\sigma^\vee}}$, we have a
surjective morphism of graded $S_n$-algebras:
$$\mathbb{C}[\overline{G \cdot M_{\sigma^\vee}} \cap \mathfrak{h} ]
\rightarrow \mathbb{C}[(N_G(T)\cdot \mathcal N_{\sigma^\vee}) \cap
\mathfrak{h}].$$ 
Thus, we get  a surjective map of graded $S_n$-algebras:
$$\phi : {\rm H}^*(X_\sigma) \rightarrow \mathbb{C}[(N_G(T)\cdot \mathcal
N_{\sigma^\vee}) \cap \mathfrak{h}].$$ 
To prove the theorem, in view of  Proposition
\ref{universal property}, it suffices to show that the top
degree $d$ of the graded algebra $\mathbb{C}[(N_G(T)\cdot \mathcal
N_{\sigma^\vee}) \cap \mathfrak{h}]$ is $d_\sigma$. As $\phi$ is
surjective, we already know that $d \leq d_\sigma$. For the other
inequality, consider the graded surjective map
$$\mathbb{C}[(N_G(T)\cdot \mathcal
N_{\sigma^\vee}) \cap \mathfrak{h}] \rightarrow \mathbb{C}
[\mathcal N_{\sigma^\vee} \cap \mathfrak{h}],$$ 
obtained by the $\mathbb{C}^*$-equivariant embedding 
$ \mathcal
N_{\sigma^\vee} \subset N_G(T)\cdot \mathcal
N_{\sigma^\vee} $.  
 By a
classical result by Kostant \cite{Kostant},  $ \mathbb{C}
[\mathcal N_{\sigma^\vee} \cap \mathfrak{h}]$ is isomorphic with the
 cohomology ${\rm
H}^*(L_{\sigma^\vee}/B_{\sigma^\vee}, \mathbb{C})$ as graded algebras
(in fact, also as modules for  the Weyl group
$N_{L_{\sigma^\vee}}(T)/T$). But the top degree of ${\rm
H}^*(L_{\sigma^\vee}/B_{\sigma^\vee}, \mathbb{C})$ coincides with
the complex dimension of $L_{\sigma^\vee}/B_{\sigma^\vee}$, which is easily
seen to be equal to $d_2$ (use formula (\ref{formula})). All
together, this implies that $d \geq d_2$. But,   by Lemma
\ref{socle}, $d_2=d_\sigma$, which  ends the proof.
\end{proof}

\end{document}